\documentclass[times,10pt,twocolumn]{article}
\usepackage{latex8}
\usepackage{times}
\usepackage{amsmath,amsthm,amsfonts,amssymb,amscd}
\usepackage{fancyhdr}
\input epsf

\newtheorem{lemma}{Lemma}
\newtheorem{theorem}[lemma]{Theorem}
\newtheorem{corollary}[lemma]{Corollary}
\newtheorem{definition}{Definition}
\newtheorem{proposition}[lemma]{Proposition}
\newtheorem{question}{Question}
\numberwithin{lemma}{section}
\numberwithin{definition}{section}
\numberwithin{question}{section}
\newcommand\C{{\mathbb{C}}}
\newcommand\R{{\mathbb{R}}}
\newcommand\Quat{{\mathbb{H}}}

\newcommand\F{{\mathbb{F}}}
\newcommand\Z{{\mathbb{Z}}}
\newcommand\Tr{{\mathop\textup{Tr }}}

\begin{document}

\title{A Group-theoretic Approach to Fast Matrix Multiplication}

\author{Henry Cohn\\
Microsoft Research\\
One Microsoft Way\\
Redmond, WA 98052-6399\\
cohn@microsoft.com\\
\and
Christopher Umans\\
Department of Computer Science\\
California Institute of Technology\\
Pasadena, CA 91125\\
umans@cs.caltech.edu\\
}

\maketitle

\begin{abstract}
We develop a new, group-theoretic approach to bounding the
exponent of matrix multiplication. There are two components to
this approach: (1) identifying groups $G$ that admit a certain
type of embedding of matrix multiplication into the group algebra
$\C[G]$, and (2) controlling the dimensions of the irreducible
representations of such groups. We present machinery and examples
to support (1), including a proof that certain families of groups
of order $n^{2 + o(1)}$ support $n \times n$ matrix
multiplication, a necessary condition for the approach to yield
exponent $2$.  Although we cannot yet completely achieve both (1)
and (2), we hope that it may be possible, and we suggest potential
routes to that result using the constructions in this paper.
\end{abstract}

\Section{Introduction}

\thispagestyle{fancy}

\fancyfoot[LE,LO]{\ \\ \parbox{6.875in}{\scriptsize{\ \\ Copyright
\copyright\ 2003 IEEE. Reprinted from Proceedings of the 44th
Annual Symposium on Foundations of Computer Science. This
material is posted here with permission of the IEEE.  Such
permission of the IEEE does not in any way imply IEEE endorsement
of any of Cornell University's products or services.  Internal or
personal use of this material is permitted.  However, permission
to reprint/republish this material for advertising or promotional
purposes or for creating new collective works for resale or
redistribution must be obtained from the IEEE by writing to
pubs-permissions@ieee.org.  By choosing to view this document,
you agree to all provisions of the copyright laws protecting
it.}}}

Strassen \cite{S} made the startling discovery that one can
multiply two $n \times n$ matrices in only $O(n^{2.81})$ field
operations, compared with $2n^3$ for the standard algorithm.  This
immediately raises the question of the exponent of matrix
multiplication: what is the smallest number $\omega$ such that for
each $\varepsilon>0$, matrix multiplication can be carried out in
at most $O(n^{\omega+\varepsilon})$ operations? Clearly $\omega
\ge 2$.  It is widely believed that $\omega=2$, but the best bound
known is $\omega < 2.38$, due to Coppersmith and Winograd
\cite{CW}, following a sequence of improvements to Strassen's
original algorithm (see \cite[p.~420]{BCS} for the history). It
is known that all the standard linear algebra problems (for
example, computing determinants, solving systems of equations,
inverting matrices, computing LUP decompositions---see Chapter~16
of \cite{BCS}) have the same exponent as matrix multiplication,
which makes $\omega$ a fundamental number for understanding
algorithmic linear algebra. In addition, there are non-algebraic
algorithms whose complexity is expressed in terms of $\omega$
(see, e.g., Section~16.9 in \cite{BCS}).

Several fairly elaborate techniques for bounding $\omega$ are
known, but since 1990 nobody has been able to improve on them. In
this paper:
\begin{itemize}
\item We develop a new approach to bounding $\omega$ that imports the
problem into the domain of group theory and representation theory.
The approach is relatively simple and almost entirely separate
{}from the existing machinery built up since Strassen's original
algorithm.

\item We demonstrate the feasibility of the group theory aspect of the
approach by identifying a family of groups for which a parameter
that mirrors $\omega$ approaches 2. We also exhibit techniques
for bounding this critical parameter and prove non-trivial bounds
for a number of diverse groups and group families.

\item We pose a question in representation theory (Question~\ref{fundamentalq}
below) that represents a potential barrier to directly obtaining
non-trivial bounds on $\omega$ using this approach. We do not
know the answer to this question. A positive answer would
illuminate a path that might lead to $\omega = 2$ using the
techniques that we present in this paper.
\end{itemize}

Our approach is reminiscent of a question asked by Coppersmith
and Winograd (in Section~11 of \cite{CW}) about avoiding ``three
disjoint equivoluminous subsets'' in abelian groups, which would
lead to $\omega=2$ if it has a positive answer. However, our
technique is completely different, and our framework seems to
have more algebraic structure to make use of (whereas theirs is
more combinatorial).

\SubSection{Analogy with fast polynomial multiplication}

There is a close analogy between the framework we propose in this
paper and the well-known algorithm for multiplying two degree $n$
polynomials in $O(n\log n)$ operations using the Fast Fourier
Transform (FFT). In this section we elucidate this analogy to give
a high-level description of our technique.

Suppose we wish to multiply the polynomials $A(x) = \sum_{i =
0}^{n-1} a_i x^i$ and $B(x) = \sum_{i = 0}^{n-1} b_i x^i$. The
naive way to do this is to compute $n^2$ products of the form
$a_ib_j$, and from these the $2n-1$ coefficients of the product
polynomial $A(x)\cdot B(x)$. Of course a far better algorithm is
possible; we describe it below in language that easily translates
into our framework for matrix multiplication.

\fancyfoot[CE,CO]{\thepage}
\fancyfoot[LE,LO]{}

Let $G$ be a group and let $\C[G]$ be the group algebra---that is,
every element of $\C[G]$ is a formal sum $\sum_{g \in G} a_g g$
with $a_g \in \C$, and the product of two such elements is
$$
\left(\sum_{g \in G}a_g g \right) \cdot \left ( \vphantom{\sum_{g
\in G}a_g g} \sum_{h \in G}b_h h \right ) = \sum_{f \in G}\left
(\sum_{gh = f} a_gb_{h} \right ) f.
$$

We often identify the element $\sum_{g \in G} a_g g$ with the
vector of its coefficients. If $G$ is the cyclic group of order
$m$, then the product of two elements $a = (a_g)_{g \in G}$ and
$b = (b_g)_{g \in G}$ is a {\em cyclic convolution} of the
vectors $a$ and $b$. The important observation is that a cyclic
convolution is almost what is needed to compute the coefficients
of the product polynomial $A(x)\cdot B(x)$---the only problem is
that it wraps around.  To avoid this problem, we embed $A(x)$ and
$B(x)$ as elements $\bar{A}, \bar{B} \in \C[G]$ as follows: Let
$z$ be a generator of $G$, which we assume to be a cyclic group
of order $m > 2n-1$, and define
\[
\bar{A} = \sum_{i = 0}^{n-1}a_iz^i \qquad \textup{and}\qquad
\bar{B} = \sum_{i = 0}^{n-1}b_iz^i.
\]
Since the group size $m$ is large enough to avoid wrapping
around, we can read off the coefficients of the product
polynomial from the element $\bar{A}\bar{B} \in \C[G]$: the
coefficient of $x^i$ in $A(x)B(x)$ is the coefficient of the
group element $z^i$ in $\bar{A}\bar{B}$. This is a wordy account
of a so-far simple correspondence, but the payoff is near. The
{\em Discrete Fourier Transform} (DFT) for $\C[G]$ is an
invertible linear transformation $D:\C[G] \rightarrow \C^{|G|}$,
which turns multiplication in $\C[G]$ into pointwise
multiplication of vectors in $\C^{|G|}$. We can therefore compute
the product $\bar{A}\bar{B}$ by first computing $D(\bar{A})$ and
$D(\bar{B})$ and then computing the inverse DFT of their
pointwise product. Thus, using the $O(m \log m)$ Fast Fourier
Transform algorithm, we can perform multiplication in $\C[G]$
(and therefore polynomial multiplication, via the embedding
above) in $O(m \log m)$ operations.

One of the main results of the present paper is that {\em matrix
multiplication can be embedded into group algebra multiplication
in an analogous way}. The embedding is not as simple as the
embedding of polynomial multiplication, but it has a natural and
clean description in terms of a property of subsets of $G$ (which
we often take to be subgroups). In particular, if $S, T$, and $U$
are subsets of $G$ and $A = (a_{s, t})_{s \in S, t \in T}$ and $B
= (b_{t, u})_{t \in T, u \in U}$ are $|S| \times |T|$ and $|T|
\times |U|$ matrices, respectively, then we define
\[\bar{A} = \sum a_{s, t}s^{-1}t \qquad\textup{and}\qquad \bar{B}
= \sum b_{t, u} t^{-1}u.\] If $S, T, U$ satisfy the {\em triple
product property} (see Definition~\ref{definition:realize}), then
we can read off the entries of the product matrix $AB$ from
$\bar{A}\bar{B} \in \C[G]$: entry $(AB)_{s, u}$ is simply the
coefficient of the group element $s^{-1}u$.

In the case of polynomial multiplication, the simplicity of the
embedding obscures the fact that if $G$ is too large (e.g., if
$|G| = n^2$ rather than $O(n)$), then the benefit of the entire
scheme is destroyed. Avoiding this pitfall turns out to be the
main challenge in the new setting. We wish to embed matrix
multiplication into a group algebra over a {\em small} group $G$,
as the size of $G$ is a lower bound on the complexity of
multiplication in $\C[G]$. It is not surprising, for example,
that $n \times n$ matrix multiplication can be embedded into the
group algebra of a group of order $n^3$. We show that abelian
groups cannot beat $n^3$ and {\em we identify families of
non-abelian groups of size $n^{2 + o(1)}$ that admit such an
embedding.}

It might seem that this result together with the above trick for
performing group algebra multiplication (i.e., taking the DFT,
multiplying in the Fourier domain, and transforming back) would
imply that $\omega = 2$. There are, however, two complications
introduced by the fact that we are forced to work with non-abelian
groups. The first is that we know of fast algorithms to compute
the DFT only for limited classes of non-abelian groups (see
Section~13.5 in \cite{BCS}). However, the DFT is linear, and
because of the recursive structure of divide and conquer matrix
multiplication algorithms, linear transformations applied before
and after the recursive step are ``free.'' For example, in
Strassen's original matrix multiplication algorithm, the number
of matrix additions and scalar multiplications in the recursive
step does not affect the bound on $\omega$. So this potential
complication is in fact no problem at all.

The second complication is that for $\C[G]$ when $G$ is
non-abelian, multiplication in the Fourier domain is {\em not}
simply pointwise multiplication of vectors in $\C^{|G|}$. Instead
it is {\em block-diagonal matrix multiplication}, where the
dimensions of the blocks are the dimensions of the irreducible
representations of $G$. We thus obtain a reduction of $n \times n$
matrix multiplication to a number of smaller matrix
multiplications of varying sizes, which gives rise to an
inequality involving the exponent $\omega$ of matrix
multiplication. If the size of $G$ were exactly $n^2$, then this
inequality would imply that $\omega = 2$. However, the smallest
one can make $|G|$ is $n^{2 + o(1)}$, and then the question of
whether the inequality implies $\omega = 2$ turns on the
representation theory of $G$. We show that when $|G| = n^{2 +
o(1)}$, even slight control over the dimension of the largest
irreducible representation is sufficient to achieve $\omega = 2$.
Some control is necessary to avoid trivialities such as reducing
to an even larger matrix multiplication problem. We can achieve
that much control; the issue of whether it is possible to achieve
more control is the subject of Question~\ref{fundamentalq}.

\SubSection{Outline}

Following some preliminaries below,
Sections~\ref{section:realizing} through~\ref{section:bounds} are
devoted to outlining our approach. In
Sections~\ref{section:linear} and~\ref{variety}, we show that a
variety of different types of groups support matrix multiplication
within our framework, and in the process demonstrate a number of
useful proof techniques. Section~\ref{section:linear} highlights
linear groups, whose representation theory makes them especially
attractive for our purposes. Section~\ref{section:Lie} describes
a parallel with Lie groups and gives a construction that suggests
that finite linear groups may indeed be a fruitful line of
inquiry. In Section~\ref{section:wreath} we consider wreath
product constructions, and in Section~\ref{section:direct} we use
the combinatorial notion of Sperner capacity to demonstrate the
surprising fact that the $k$-fold direct product of a group may
support $n^k \times n^k$ matrix multiplication even when the group
itself fails to support $n \times n$ matrix multiplication. This
suggests a potential route to answering
Question~\ref{fundamentalq} in the affirmative. We end by
mentioning some open problems and variants of our overall
approach in Section~\ref{section:conclusions}.

\SubSection{Preliminaries}

Let $\langle n, m, p \rangle$ denote the structural tensor for
rectangular matrix multiplication of $n \times m$ by $m \times p$
matrices, and let $R$ denote the tensor rank function. (See
\cite{BCS} for background on matrix multiplication and tensors.
We will use this material only in the proof of
Theorem~\ref{theorem:bound}.)  We will typically work over the
field of complex numbers; if we use another field $F$, we will
write $\langle n, m, p \rangle_F$. As usual $\omega$ will denote
the exponent of matrix multiplication over $\C$.

We will use the following basic fact from representation theory:
the group algebra $\C[G]$ of a finite group $G$ decomposes as the
direct product
$$
\C[G] \cong  \C^{d_1 \times d_1} \times \dots \times \C^{d_k
\times d_k}
$$
of matrix algebras of orders $d_1,\dots,d_k$.  These numbers are
called the character degrees of $G$, or the dimensions of the
irreducible representations.  It follows from computing the
dimensions of both sides that $|G| = \sum_i d_i^2$.  See
\cite{JL} and \cite{H} for background on representation theory.

\Section{Realizing matrix multiplication via groups}
\label{section:realizing}

In this section we describe the embedding of matrix multiplication
into group algebra multiplication, and we identify a property of
groups $G$ that implies that the group algebra of $G$ admits such
an embedding. If $S$ is a subset of a group, let $Q(S)$ denote
the right quotient set of $S$, i.e.,
$$
Q(S) = \{s_1 s_2^{-1} : s_1,s_2 \in S\}.
$$

\begin{definition}
\label{definition:realize} A group $G$ {\em realizes} $\langle
n_1,n_2,n_3 \rangle$ if there are subsets $S_1,S_2,S_3 \subseteq
G$ such that $|S_i| = n_i$, and for $q_i \in Q(S_i)$, if
$$
q_1q_2q_3 = 1
$$
then $q_1=q_2=q_3=1$. We call this condition on $S_1,S_2,S_3$ the
{\em triple product property}. If we wish to emphasize the
specific subsets, we say that $G$ {\em realizes $\langle
n_1,n_2,n_3 \rangle$ through} $S_1,S_2,S_3$.
\end{definition}

In most of our examples, matrix multiplication will be realized
through subgroups $H_1$, $H_2$, $H_3$ of $G$, rather than
arbitrary subsets.  In that case, the triple product property is
especially simple, because $Q(H_i) = H_i$: it states that if
$h_1h_2h_3=1$ with $h_i \in H_i$, then $h_1=h_2=h_3=1$.  An
equivalent formulation replaces $h_1h_2h_3=1$ with $h_1h_2=h_3$.

Perhaps the simplest example comes from the product $C_{n} \times
C_m \times C_p$ of cyclic groups, which clearly realizes $\langle
n,m,p \rangle$ through $C_n \times \{1\} \times \{1\}$, $\{1\}
\times C_m \times \{1\}$, and $\{1\} \times \{1\} \times C_p$.  We
will see a number of less trivial examples shortly.

\begin{lemma}
\label{lemma:permute} If $G$ realizes $\langle
n_1,n_2,n_3\rangle$, then it does so for every permutation of
$n_1,n_2,n_3$.
\end{lemma}

\begin{proof}
Suppose $G$ realizes $\langle n_1,n_2,n_3\rangle$ through
$S_1,S_2,S_3$, and suppose $s_i,s_i' \in S_i$. We need to show
that the order in which $1$, $2$, and $3$ appear in the equation
$$
s_1's_1^{-1} s_2's_2^{-1} s_3's_3^{-1} = 1
$$
is irrelevant.  Conjugating by $s_1' s_1^{-1}$ shows that it is
equivalent to
$$
s_2's_2^{-1} s_3's_3^{-1} s_1's_1^{-1} = 1,
$$
so we can perform a cyclic shift.  To get a transposition, we
take the inverse of the initial equation, which yields
$$
s_3s_3'^{-1} s_2s_2'^{-1} s_1s_1'^{-1} = 1,
$$
i.e., a transposition of $1$ with $3$ (the roles of $s$ and $s'$
have been reversed, but that is irrelevant).  These two
permutations generate all permutations of $\{1,2,3\}$.
\end{proof}

\begin{lemma}
\label{lemma:shortexact} If $N$ is a normal subgroup of $G$ that
realizes $\langle n_1, n_2, n_3\rangle$ and $G/N$ realizes
$\langle m_1, m_2, m_3\rangle$, then $G$ realizes $\langle
n_1m_1, n_2m_2, n_3m_3\rangle$.
\end{lemma}

\begin{proof}
Suppose $N$ realizes $\langle n_1, n_2, n_3\rangle$ through
$S_1,S_2,S_3$, and suppose $T_1,T_2,T_3$ are lifts to $G$ of the
three subsets of $G/N$ that realize $\langle m_1, m_2,
m_3\rangle$. Then we claim that $G$ realizes $\langle n_1m_1,
n_2m_2, n_3m_3\rangle$ through the pointwise products
$S_1T_1,S_2T_2,S_3T_3$.  We need to check that for $s_i,s_i' \in
S_i$ and $t_i,t_i' \in T_i$,
$$
(s_1't_1')(s_1t_1)^{-1} (s_2't_2')(s_2t_2)^{-1}
(s_3't_3')(s_3t_3)^{-1} = 1
$$
iff $s_i=s_i'$ and $t_i=t_i'$ for all $i$.  If we reduce this
equation modulo $N$, we find that $t_i=t_i'$ modulo $N$, and
hence also in $G$.  The equation in $G$ then becomes
$$
s_1's_1^{-1} s_2's_2^{-1} s_3's_3^{-1} = 1,
$$
{}from which we deduce $s_i=s_i'$, as desired.
\end{proof}

One useful special case of Lemma~\ref{lemma:shortexact} is that if
$G_1$ realizes $\langle n_1, m_1, p_1\rangle$ and $G_2$ realizes
$\langle n_2, m_2, p_2\rangle$, then $G_1 \times G_2$ realizes
$\langle n_1n_2, m_1m_2, p_1p_2\rangle$.

Our first theorem describes the embedding of matrix multiplication
into group algebra multiplication:

\begin{theorem}
\label{theorem:reduction} Let $F$ be any field.  If $G$ realizes
$\langle n,m,p \rangle$, then the number of field operations
required to multiply $n \times m$ with $m \times p$ matrices over
$F$ is at most the number of operations required to multiply two
elements of $F[G]$.  Furthermore, $\langle n,m,p\rangle_F \le
F[G]$.
\end{theorem}

For the definition of the restriction relation $\le$ in the last
sentence, see Section~14.3 of \cite{BCS}.

\begin{proof}
Let $G$ realize $\langle n,m,p \rangle$ through subsets $S,T,U$.
Suppose $A$ is an $n \times m$ matrix, and $B$ is an $m \times p$
matrix.  We will index the rows and columns of $A$ with the sets
$S$ and $T$, respectively, those of $B$ with $T$ and $U$, and
those of $AB$ with $S$ and $U$.

Consider the product
$$
\left(\sum_{s\in S, t \in T} A_{st} s^{-1} t\right)
\left(\sum_{t' \in T, u \in U} B_{t'u} t'^{-1} u\right)
$$
in the group algebra.  We have
$$
(s^{-1} t) (t'^{-1} u) = s'^{-1} u'
$$
iff $s=s'$, $t=t'$, and $u=u'$, so the coefficient of $s^{-1} u$
in the product is
$$
\sum_{t \in T} A_{st} B_{tu} = (AB)_{su}.
$$
Thus, one can simply read off the matrix product from the group
algebra product by looking at the coefficients of $s^{-1} u$ with
$s \in S, u \in U$, and the assertions in the theorem statement
follow.
\end{proof}

\Section{The pseudo-exponent}

The {\em pseudo-exponent} of a group measures the quality of the
embedding afforded by Theorem \ref{theorem:reduction} in a single,
well-behaved parameter, which in some ways mirrors the exponent
$\omega$ of matrix multiplication.

\begin{definition}
The pseudo-exponent $\alpha(G)$ of a non-trivial finite group $G$
is the minimum of
$$
\frac{3 \log |G|}{\log nmp}
$$
over all $n,m,p$ (not all $1$) such that $G$ realizes $\langle
n,m,p \rangle$.  The pseudo-exponent of the trivial group is $3$.
\end{definition}

When it is clear from the context which group is intended, we
often write $\alpha$ instead of $\alpha(G)$.  Note that in the
special case that $G$ realizes $\langle n, n, n \rangle$, its
pseudo-exponent satisfies $\alpha \le \log_n |G|$.  In general,
if $G$ realizes $\langle n,m,p \rangle$, then
$$
\alpha \le \log_{\sqrt[3]{nmp}} |G|.
$$

\begin{lemma}
The pseudo-exponent of a finite group $G$ is always greater than
$2$ and at most $3$. If $G$ is abelian, then it is exactly $3$.
\end{lemma}

\begin{proof}
The upper bound of $3$ is trivial: use the subgroups $H_1=H_2 =
\{1\}$ and $H_3=G$.

For the lower bounds, suppose $G$ realizes $\langle n_1,n_2,n_3
\rangle$ (with $n_1n_2n_3>1$) through subsets $S_1,S_2,S_3$.  It
follows from the definition of realization that the map $(x,y)
\mapsto x^{-1} y$ is injective on $S_1 \times S_2$ and its image
intersects the quotient set $Q(S_3)$ only in the identity.  Thus,
$|G| \ge n_1n_2$, and $|G| > n_1n_2$ unless $n_3=1$.  Similarly,
$|G| \ge n_2n_3$ with equality only if $n_1=1$, and $|G| \ge
n_1n_3$ with equality only if $n_2=1$. Thus, $|G|^3 >
(n_1n_2n_3)^2$, so $\alpha(G) > 2$.

If $G$ is abelian, then the product map $S_1 \times S_2 \times S_3
\to G$ must be injective, so $|G| \ge n_1n_2n_3$ and $\alpha(G)
\ge 3$.
\end{proof}

The pseudo-exponent is well-behaved with respect to group
extensions:

\begin{lemma}
\label{lemma:upperbound} If $N$ is a normal subgroup of $G$, then
$\alpha(G) \le \max(\alpha(N),\alpha(G/N))$.
\end{lemma}

\begin{proof}
Suppose $N$ realizes $\langle n_1, n_2, n_3\rangle$ and $G/N$
realizes $\langle m_1, m_2, m_3\rangle$.  Then
Lemma~\ref{lemma:shortexact} implies that the pseudo-exponent of
$G$ is at most
$$
\frac{3\log |G|}{\log n_1m_1n_2m_2n_3m_3} = \frac{3\log |N| +
3\log |G/N| }{\log n_1n_2n_3 + \log m_1m_2m_3},
$$
which is bounded above by the larger of
$$
\frac{3\log |N|}{\log n_1n_2n_3} \qquad \textup{and}\qquad
\frac{3\log |G/N|}{\log m_1m_2m_3},
$$
as desired.
\end{proof}

Non-abelian groups can have pseudo-exponent less than $3$.  The
smallest example is the symmetric group $S_3$ on $3$ elements. It
realizes $\langle 2, 2, 2 \rangle$ through its three subgroups of
order $2$, so it has pseudo-exponent at most $\log_2 6$ (and one
can check that it is exactly $\log_2 6$). Next, we generalize
this construction to show that it is possible to come arbitrarily
close to pseudo-exponent $2$, as follows.

\begin{figure}
\begin{center}
\leavevmode \epsfbox[-2 -2 82 92]{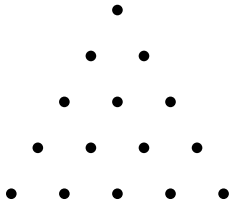}
\end{center}
\caption{A triangular array of points.} \label{fig-triangle}
\end{figure}

Given a triangular array of points in the plane, as in
Figure~\ref{fig-triangle}, we consider the group of permutations
of the points, together with three subgroups, one for each side
of the triangle.  Each subgroup permutes the set of points on
each line parallel to its side of the triangle.  The proof of
Theorem~\ref{theorem:pseudoexponent2}, while not phrased in
geometric terms, shows that these subgroups satisfy the triple
product property.

\begin{theorem}
\label{theorem:pseudoexponent2} The pseudo-exponent of
$S_{n(n+1)/2}$ is at most
$$
2 +\frac{2-\log 2}{\log n} + O\left(\frac{1}{(\log n)^2}\right).
$$
\end{theorem}

\begin{proof}
There are $n(n+1)/2$ triples $(a,b,c)$ with $a,b,c \ge 0$ and
$a+b+c=n-1$.  We view $S_{n(n+1)/2}$ as the group of permutations
of these triples.  Let $H_i$ be the subgroup that fixes the $i$-th
coordinate.  The size of this subgroup is $1!2!\dots n!$, so the
pseudo-exponent bound is
$$
\frac{\log (n(n+1)/2)!}{\log 1!2!\dots n!} = 2 + \frac{2-\log
2}{\log n} + O\left(\frac{1}{(\log n)^2}\right),
$$
assuming these subgroups satisfy the triple product property.  For
that, we need to prove that if $h_1h_2h_3=1$ with $h_i \in H_i$,
then $h_1=h_2=h_3=1$.

Suppose $h_1h_2h_3=1$ with $h_i \in H_i$.   We will order the
triples lexicographically, so that $(0,0,n-1)$ is the smallest
triple and $(n-1,0,0)$ is the largest, and prove by induction
using this ordering that $h_1$, $h_2$, and $h_3$ fix every triple.

Suppose all triples smaller than $(a,b,c)$ are fixed by each of
$h_1, h_2, h_3$ (in the base case, the set of such triples is
empty). The permutation $h_3$ cannot send $(a,b,c)$ to a smaller
triple, since all smaller triples are fixed points, so $h_3$ must
send it to $(a+i,b-i,c)$ with $i \ge 0$.  Then $h_2$ sends that
to $(a+i+j,b-i,c-j)$ for some $j$. The only way $h_1$ can return
to $(a,b,c)$ is if $i+j=0$, so that must be the case. However,
$h_1$ fixes $(a,b-i,c+i)$ for $i>0$ (since such a triple is
smaller than $(a,b,c)$), so we must have $i=0$.  It follows that
$(a,b,c)$ is fixed by each of $h_1,h_2,h_3$, so by induction all
triples are fixed and hence $h_1=h_2=h_3=1$.
\end{proof}

The same holds for all symmetric groups, since one can look at the
largest subgroup of the form $S_{n(n+1)/2}$.

\Section{Relating the pseudo-exponent to $\omega$}
\label{section:bounds}

In this section we relate the pseudo-exponent $\alpha$ to the
exponent of matrix multiplication $\omega$.  As with many of the
results since Strassen's algorithm, our main theorems are stated
as bounds on $\omega$, rather than explicit algorithms, but of
course algorithms are implicit in the proofs.

\begin{theorem}
\label{theorem:bound} Suppose $G$ has pseudo-exponent $\alpha$,
and the character degrees of $G$ are $\{d_i\}$. Then
$$
|G|^{\omega/\alpha} \le \sum_{i} d_i^{\omega}.
$$
\end{theorem}

The intuition is simple: the problem of multiplying matrices of
size $|G|^{1/\alpha}$ reduces to multiplication in $\C[G]$, which
is equivalent to multiplying a collection of matrices of sizes
$d_i$.  These multiplications should take about $d_i^\omega$
operations, so $\sum_i d_i^\omega$ should be an approximate upper
bound for the number of operations required to multiply matrices
of size $|G|^{1/\alpha}$, i.e., roughly $|G|^{\omega/\alpha}$. It
is convenient that when one makes this idea precise, these crude
approximations become exact bounds.

\begin{proof}
Suppose $G$ realizes $\langle n, m, p \rangle$ with $nmp =
|G|^{3/\alpha}$ (it follows from the definition of the
pseudo-exponent that $G$ realizes such a tensor). By
Theorem~\ref{theorem:reduction},
\begin{equation}
\label{equation:reduction} \langle n, m, p \rangle \le \C[G]
\simeq \bigoplus_i \langle d_i, d_i, d_i \rangle.
\end{equation}
We will need two facts about the rank of matrix multiplication:
for all $n',m',p'$,
$$
(n'm'p')^{\omega/3} \le R(\langle n', m', p' \rangle)
$$
(Proposition~15.5 in \cite{BCS}), and for each $\varepsilon>0$
there exists $C>0$ such that for all $k$,
$$
R(\langle k, k, k \rangle) \le C k^{\omega+\varepsilon}
$$
(Proposition~15.1 in \cite{BCS}).

The $\ell$-th tensor power of \eqref{equation:reduction} is
$$
\langle n^\ell, m^\ell, p^\ell \rangle \le
\bigoplus_{i_1,\dots,i_\ell} \langle d_{i_1}\dots
d_{i_\ell},d_{i_1}\dots d_{i_\ell}, d_{i_1}\dots d_{i_\ell}
\rangle,
$$
if we use
$$
\langle n_1,m_1,p_1\rangle \otimes \langle n_2,m_2,p_2\rangle
\simeq \langle n_1n_2,m_1m_2,p_1p_2 \rangle.
$$
It follows from taking the rank of both sides that
$$
|G|^{\ell\omega/\alpha} \le C \left(\sum_i
d_i^{\omega+\varepsilon}\right)^\ell,
$$
and if we take the $\ell$-th root and let $\ell$ go to infinity,
then we deduce that
$$
|G|^{\omega/\alpha} \le \sum_{i} d_i^{\omega+\varepsilon}.
$$
Finally, because this inequality holds for all $\varepsilon>0$, it
must hold for $\varepsilon=0$ as well, by continuity.
\end{proof}

Notice that if $\alpha(G)$ were $2$, then this theorem would imply
that $\omega=2$ (using $\sum_i d_i^2 = |G|$, the Cauchy-Schwarz
inequality, and the fact that every non-trivial group has at least
two irreducible representations). In general, though, we need to
control the character degrees of $G$. The maximum possible
character degree for any non-trivial group is $(|G|-1)^{1/2}$; we
show below that an upper bound of $|G|^{1/2 - \varepsilon}$ for
fixed $\varepsilon > 0$ would be sufficient to obtain $\omega =
2$ from a family of groups with pseudo-exponent approaching $2$
(and that even a much weaker bound suffices).

We define $\gamma(G)$, or simply $\gamma$ when $G$ is clear from
the context, so that $|G|^{1/\gamma}$ is the maximum character
degree of $G$ ($\gamma(G)=\infty$ if $G$ is abelian). Ideally,
we'd like the exponent of matrix multiplication $\omega$ to be
bounded above by the pseudo-exponent $\alpha$. The following
corollary shows that in the region near 2, this actually happens,
with a correction factor that depends on $\gamma$.

\begin{corollary}
\label{corollary:upperbound} Let $G$ be a finite group. If
$\alpha(G) < \gamma(G)$, then
$$
\omega \le \alpha \left (\frac{\gamma - 2}{\gamma - \alpha} \right
).
$$
\end{corollary}

\begin{proof}
Let $\{d_i\}$ denote the character degrees.  Then by
Theorem~\ref{theorem:bound},
\begin{eqnarray*}
|G|^{\omega/\alpha} &\le& \sum_i d_i^{\omega -2}d_i^2\\
&\le& |G|^{(\omega - 2)/\gamma}\sum_i d_i^2\\
& = & |G|^{1+(\omega - 2)/\gamma},
\end{eqnarray*}
which implies $\omega(1/\alpha-1/\gamma) \le 1-2/\gamma.$
Dividing by $1/\alpha - 1/\gamma$ (which is positive by
assumption) yields the stated result.
\end{proof}

Like $\alpha(G)$, we have $\gamma(G)>2$ for all $G$, and
Corollary~\ref{corollary:upperbound} shows that our approach
amounts to a race between $\alpha(G)$ and $\gamma(G)$ to see
which approaches $2$ faster. The most attractive form of this
corollary is the following special case:

\begin{corollary}
\label{cor:race} Suppose there exists a family $G_1,G_2,\dots$ of
finite groups such that $\alpha(G_i) = 2+o(1)$ as $i \to \infty$,
and furthermore $\alpha(G_i)-2 = o(\gamma(G_i)-2)$.  Then the
exponent of matrix multiplication is $2$.
\end{corollary}

These corollaries are weakenings of Theorem~\ref{theorem:bound},
the advantage being that they only require knowledge of
$\gamma(G)$, which is typically easier to work with than the
complete set of character degrees that is required for
Theorem~\ref{theorem:bound}.

It is reasonable to ask whether the requirement $\alpha < \gamma$
which occurs in Corollary~\ref{corollary:upperbound} is
necessary. It turns out that it is, because if $\alpha \ge
\gamma$, then for all $\omega>0,$
$$
|G|^{\omega/\alpha} \le |G|^{\omega/\gamma} \le \sum_i d_i^\omega,
$$
where the second inequality holds because $|G|^{1/\gamma} = d_i$
for some $i$. Then the inequality in Theorem~\ref{theorem:bound}
holds even for $\omega = 3$. The necessity of $\alpha < \gamma$
makes perfect sense, because when it fails to hold, the approach
amounts to a reduction of matrix multiplication to several
instances, one of which is as large as the original instance. In
fact, the construction in the proof of
Theorem~\ref{theorem:pseudoexponent2} succumbs to this problem:
there we proved that $\alpha(S_{n(n+1)/2}) \le 2+O(1/\log n)$,
but it turns out that $\gamma(S_{n(n+1)/2}) = 2 + \Theta(1/(n\log
n))$ (see \cite{VK}).  However, there exist non-abelian groups
for which $\alpha < \gamma$ and $\alpha < 3$; one example is the
group in Proposition~\ref{proposition:order80} below.

If we {\em do} have access to the complete set of character
degrees then there is a relatively simple condition to check to
determine whether the inequality in Theorem~\ref{theorem:bound}
yields a non-trivial bound on $\omega$. The condition is that
$|G|^{3/\alpha} > \sum_i{d_i^3}$. To see this observe that the
inequality in Theorem~\ref{theorem:bound} is equivalent to
\begin{equation} \label{eqn:log}
\frac{\omega}{\alpha} \log |G| \le \log \sum_i d_i^\omega.
\end{equation}
The right-hand side is convex as a function of $\omega$, and the
left-hand side is linear. Furthermore, as $\omega \to \infty$,
the right-hand side is asymptotic to
$$
\frac{\omega}{\gamma} \log |G|,
$$
which is smaller than the left-hand side when $\alpha < \gamma$
(which is the non-trivial case).  Therefore \eqref{eqn:log} gives
no information about $\omega$ in the interval $[2,3]$ unless it
rules out $\omega=3$, which is equivalent to the above stated
condition. We do not have examples of groups meeting this
condition.

We are thus led to pose the following question in representation
theory:

\begin{question}
\label{fundamentalq} Does there exist a finite group that realizes
$\langle n,m,p \rangle$ and has character degrees $\{d_i\}$ such
that
$$
nmp > \sum_i d_i^3?
$$
\end{question}

It is possible that there is a theorem in representation theory
that implies that the answer to this question is ``no.'' In that
case the approach we have outlined cannot be used directly to
obtain bounds on $\omega$; however even in this case there are
variants of our approach that would not be ruled out (see, e.g.,
Subsection~\ref{section:extensions}). On the other hand, a
positive answer might point the direction to a proof that $\omega
= 2$ using our approach: it would seem strange if the best bound
groups could prove were some constant strictly between~$2$
and~$3$, and the condition in Corollary~\ref{cor:race} for
$\omega=2$ feels very natural.

\Section{Linear groups} \label{section:linear}

Matrix groups over finite fields are an important class of finite
groups. They are especially attractive for our purposes because
their representation sizes, as measured by $\gamma$, are well
behaved. We will focus on the case of $SL_n(\F_q)$ for simplicity,
although we see no reason why it should perform better than other
linear groups.  If $n>1$ is held fixed, $\gamma(SL_n(\F_q))$
approaches $2+2/n$ as $q$ tends to infinity (which can be deduced
from \cite{Green}, according to a private
communication from G.\ Lusztig). Thus, if one could prove that
$\alpha(SL_n(\F_q)) = 2+o(1)$ for some fixed $n$, then
Corollary~\ref{corollary:upperbound} would imply $\omega=2$. Even
if one lets $n$ grow, one might still hope that $\alpha$ would
tend to $2$ faster than $\gamma$.  We cannot prove that $\alpha$
even approaches $2$ at all as $n,q \to \infty$, but comparison
with Theorem~\ref{theorem:slnLie} below suggests that it does. In
this section we concentrate on the case of $SL_2(\F_q)$.

For later reference, we collect here the character degrees of
$SL_2(\F_q)$:
\begin{table}[h]
\begin{tabular}{c|c|c}
Degree & Multiplicity ($q$ odd) & Multiplicity ($q$ even)\\ \hline
$q+1$ & $(q-3)/2$ & $(q-2)/2$\\
$q$ & $1$ & $1$\\
$q-1$ & $(q-1)/2$ & $q/2$\\
$(q+1)/2$ & $2$ & $0$\\
$(q-1)/2$ & $2$ & $0$\\
$1$ & $1$ & $1$
\end{tabular}
\end{table}

(See Exercise~28.2 and its solution in \cite{JL} for $q$ even, and
\cite{LR} for $q$ odd, but note that \cite{LR} has a typo in the
multiplicity for degree $q+1$ at the bottom of the first column
on page~122.)

\begin{proposition}
\label{proposition:sl2} The group $SL_2(\F_q)$ of order $q^3-q$
realizes $\langle q, q, q \rangle$.
\end{proposition}

Unfortunately, this pseudo-exponent bound tends to $3$ as $q \to
\infty$, but at least it is always strictly better than $3$.  (We
can also prove similarly that $\alpha(SL_n(\F_q)) < 3$.)

\begin{proof}
Consider the three parabolic subgroups
$$
H_1 = \left\{ \left(
\begin{array}{cc}
1 & x\\
0 & 1 \end{array} \right) : x \in \F_q\right\},
$$
$$
H_2 = \left\{ \left(
\begin{array}{cc}
1 & 0\\
y & 1 \end{array} \right) : y \in \F_q\right\},
$$
and
$$
H_3 = \left\{ \left(
\begin{array}{cc}
1+z & z\\
-z & 1-z \end{array} \right) : z \in \F_q\right\}.
$$
We need to check that for $h_i \in H_i$, if $h_1h_2=h_3$, then
$h_1=h_2=h_3=1$. To check that, we multiply to get
$$
\left(
\begin{array}{cc} 1 & x\\ 0 & 1
\end{array} \right)\left( \begin{array}{cc} 1 & 0\\ y & 1
\end{array} \right) = \left( \begin{array}{cc}
1+xy & x\\
y & 1 \end{array} \right).
$$
That can be of the form
$$
\left(
\begin{array}{cc}
1+z & z\\
-z & 1-z \end{array} \right)
$$
only if $x=y=z=0$, as desired.
\end{proof}

One might hope that $SL_n(\F_q)$ realizes
$$
\langle q^{n(n-1)/2}, q^{n(n-1)/2}, q^{n(n-1)/2}\rangle
$$
through three conjugates of the group of upper-triangular matrices
with $1$'s on the diagonal. However, that fails for $q=2$ and
$n=3$, according to calculations using the computer program GAP
(see~\cite{GAP}); furthermore, no subgroups of these orders work
for $q=2$ and $n=3$.

\begin{proposition}
\label{proposition:sl2fq2} The group $SL_2(\F_{q^2})$ of order
$q^6-q^2$ realizes $\langle q^2, q^2, q^3-q \rangle$.
\end{proposition}

\begin{proof}
Let $x \mapsto \bar x$ denote the Frobenius automorphism of
$\F_{q^2}$ over $\F_q$.  The three subgroups we will use are
$$
H_1 = \left\{ \left(
\begin{array}{cc}
1 & x\\
0 & 1 \end{array} \right) : x \in \F_{q^2}\right\},
$$
$$
H_2 = \left\{ \left(
\begin{array}{cc}
1 & 0\\
y & 1 \end{array} \right) : y \in \F_{q^2}\right\},
$$
and
\begin{eqnarray*}
H_3 &=& SU_2(\F_q)\\
&=& \left\{ \left(
\begin{array}{cc}
a & b\\
-\bar b & \bar a \end{array} \right) : a,b \in \F_{q^2}, a\bar a
+ b\bar b = 1\right\}.
\end{eqnarray*}
Note that to check that $|H_3| = q^3-q$, one just needs to count
solutions to $a\bar a + b\bar b = 1$. For a fixed $b$ with $b\bar
b \ne 1$, there are $q+1$ corresponding choices of $a$ that work;
if $b\bar b =1$, then $a=0$.  There are $(q^2-1)-(q+1)$ non-zero
choices of $b$ with $b \bar b \ne 1$ (to which we must add
$b=0$), and $q+1$ with $b\bar b = 1$.  Thus, there are
$(q^2-q-1)(q+1)+(q+1) = q^3-q$ elements of $H_3$.

As in the previous proof, checking the triple product property
amounts to checking that
$$
\left( \begin{array}{cc}
1+xy & x\\
y & 1 \end{array} \right) = \left(
\begin{array}{cc}
a & b\\
-\bar b & \bar a \end{array} \right)
$$
implies $x=y=b=0$ and $a=1$, which is a trivial calculation.
\end{proof}

Proposition~\ref{proposition:sl2fq2} proves that
$$
\liminf_{q \to \infty} \alpha(SL_2(\F_q)) \le 18/7,
$$
which is substantially better than~$3$ but still not near $2$.
Using Theorem~\ref{theorem:bound} and the character degrees of
$SL_2(\F_q)$, one can show that if
$$
\liminf_{q \to \infty} \alpha(SL_2(\F_q)) < 9/4,
$$
then Question~\ref{fundamentalq} has a positive answer.

\Section{Lie groups} \label{section:Lie}

In the category of Lie groups, one can set up a theory parallel to
that of the previous sections. We do not know how to use it to
bound the exponent of matrix multiplication (because of course
Lie groups of positive dimension are infinite). However, we have
had more luck constructing examples using Lie groups than with
finite linear groups, and this success seems a good reason to be
optimistic about matrix groups over finite fields. All examples
involving Lie groups can be skipped by a reader who cares only
about finite groups and matrix multiplication.

Recall that $Q(S)$ denotes the right quotient set of $S$.

\begin{definition}
\label{definition:Liepseudoexponent}  Let $G$ be a Lie group,
with submanifolds $M_1, M_2, M_3$ such that for $q_i \in Q(M_i)$,
if $q_1q_2q_3=1$ then $q_1=q_2=q_3=1$. We say that $G$ has {\em
Lie pseudo-exponent} at most
$$
\frac{\dim(G)}{(\dim(M_1)+\dim(M_2)+\dim(M_3))/3}.
$$
\end{definition}

We usually take the submanifolds to be Lie subgroups.  If $G$ and
the three subgroups are algebraic groups defined over a number
field, then it is natural to ask what pseudo-exponent may be
achieved when one reduces modulo a prime ideal, to get a finite
quotient group. If the triple product property still holds, then
as the finite field size tends to infinity, the pseudo-exponent
bound of this finite group approaches the Lie pseudo-exponent.
However, the triple product property may not be preserved, as we
will show after the following theorem.

\begin{theorem}
\label{theorem:slnLie} The group $SL_n(\R)$ has Lie
pseudo-exponent at most $2+2/n$.
\end{theorem}

\begin{proof}
The three subgroups are the group $U$ of upper-triangular
matrices with $1$'s on the diagonal, the group $L$ of
lower-triangular matrices with $1$'s on the diagonal, and
$SO_n(\R)$.  Each subgroup has dimension $n(n-1)/2$, and
$SL_n(\R)$ has dimension $n^2-1$, so assuming the triple product
property holds, the Lie pseudo-exponent is at most
$$
\frac{n^2-1}{n(n-1)/2} = 2+\frac{2}{n}.
$$

Let $M \in SO_n(\R)$, $A \in U$, and $B \in L$.  We wish to prove
that if $MA=B$, then $M=A=B=I$.  Let $e_1,\dots,e_n$ be the
standard basis of $\R^n$.  We will prove by induction on $i$ that
$Me_i=e_i$.  Once we know that $M=I$, it follows that $A=B$, and
thus $A=B=I$ because $U$ and $L$ are disjoint except for the
identity.  ($A=B=I$ will also follow directly from the proof that
$M=I$.)

Let $A_i$ and $B_i$ denote the $i$-th columns of $A$ and $B$, and
denote their $j$-th entries by $A_{ij}$ and $B_{ij}$.  Note that
this indexing of rows and columns is opposite to the standard
convention, but it will be more convenient in this proof.
Because $MA=B$, we have
$$
M A_i = B_i.
$$

We start with the base case $i=1$.  Since $A$ is in $U$, we have
$A_1 = e_1$.  Thus, $|B_1| = |MA_1| = |Me_1| = |e_1| = 1$, since
$M$ is an orthogonal matrix.  Because $B_{11} = 1$, the only way
$|B_1|$ can be $1$ is if $B_1=e_1$.  Thus, $Me_1 = e_1$.

Now suppose that $Me_j = e_j$ for all $j<i$.  Because $A$ is in
$U$,
$$
A_i = e_i + \sum_{j < i} A_{ij} e_j,
$$
and because $B$ is in $L$,
$$
B_i = e_i + \sum_{j > i} B_{ij} e_j.
$$
Now the induction hypothesis implies that
$$
B_i = MA_i = Me_i + \sum_{j < i} A_{ij} e_j,
$$
so
$$
Me_i = e_i +\sum_{j > i} B_{ij} e_j - \sum_{j < i} A_{ij}e_j.
$$
Since $M$ is orthogonal, $|Me_i|=|e_i|=1$.  The coefficient of
$e_i$ in $Me_i$ is already $1$, so the other coefficients must be
zero and thus $Me_i=e_i$, as desired.
\end{proof}

The same holds for $SL_n(\C)$ with $SO_n(\R)$ replaced by $SU_n$,
but not by $SO_n(\C)$: the orthogonal matrix
$$
\left(\begin{array}{ccc} 1 & \frac{-1+i}{2} & \frac{1+i}{2}\\
1 & \frac{1+i}{2} & \frac{-1+i}{2}\\
-i & 1 & 1
\end{array}\right)
$$
equals
$$
\left(\begin{array}{ccc} 1 & 0 & 0\\
1 & 1 & 0\\
-i & \frac{1-i}{2} & 1
\end{array}\right)
\left(\begin{array}{ccc} 1 & \frac{-1+i}{2} & \frac{1+i}{2}\\
0 & 1 & -1\\
0 & 0 & 1
\end{array}\right).
$$
Of course the same obstacle arises over finite fields (a sum of
non-zero squares may vanish).

\Section{Additional examples} \label{variety}

In this section we explore a variety of different types of groups,
and prove non-trivial pseudo-exponent bounds for them. We hope
that these examples (together with the ones we have already seen)
will serve as something of a tool kit for constructing a group
that might answer Question~\ref{fundamentalq}, and possibly even a
family of groups that prove $\omega = 2$.

\SubSection{Solvable groups} \label{section:solvable}

Non-abelian simple (or almost simple) groups appear to be a
fruitful source of groups with small pseudo-exponents.  However,
solvable groups also do quite well. In this section, we will
construct solvable groups that have Lie pseudo-exponent tending
to $2$, and finite solvable groups with pseudo-exponent bounds of
$2.5$ and $2.4811\dots$ (which, GAP tells us, is the best
pseudo-exponent attained using three subgroups in any group of
order up to~$100$).

Let $F$ be a field, and $\langle,\rangle$ a symmetric bilinear
form on $F^n$.  Define multiplication in
$$
G = \{(x,y,\alpha): x,y \in F^n, \alpha \in F \}
$$
via
$$
(x,y,\alpha) (u,v,\beta) = (x+u,y+v,\alpha+\beta+2\langle u, y
\rangle),
$$
and define the three subgroups
$$
H_1 = \{(x,0,0): x \in F^n\},
$$
$$
H_2 = \{(0,y,0): y \in F^n\},
$$
and
$$
H_3 = \{(z,z,\langle z,z \rangle): z \in F^n\}.
$$

\begin{proposition}
If the only element $z \in F^n$ satisfying $\langle z,z \rangle =
0$ is $z=0$, then $H_1$, $H_2$, and $H_3$ satisfy the triple
product property. \label{prop:inner}
\end{proposition}

\begin{proof}
We simply need to check that $H_3$ avoids all elements of the form
$(x,0,0)(0,y,0) = (x,y,0)$, except when $x=y=0$.  The only way
such an element can be in $H_3$, i.e., of the form $(z,z,\langle
z,z \rangle)$, is if $x=y=z$ and $\langle z,z \rangle = 0$.  That
means $z=0$ and thus $x=y=0$, as desired.
\end{proof}

When $F = \R$, the group described above is a Heisenberg group,
and we obtain the following bound:

\begin{corollary}
In the above framework, with $F = \R$, and $\langle,\rangle$ the
standard inner product, the Lie group $G$ has Lie pseudo-exponent
at most $2 + 1/n$.
\end{corollary}

\begin{proof}
It is clear that Proposition \ref{prop:inner} is satisfied; the
group dimension is $2n + 1$, and the three subgroups each have
dimension $n$.
\end{proof}

When $F$ is a finite field, the group described above is an
extraspecial group, and we obtain the following bound:

\begin{corollary}
In the above framework, with $F = \F_q$ of odd characteristic, $n
= 2$, and $\langle x, y \rangle = x_1y_1 - wx_2y_2$ for some $w
\in F$ that is not a square, the finite group $G$ has
pseudo-exponent at most $2.5$.
\end{corollary}

Here, $x_i$ denotes the $i$-th coordinate of the vector $x$.

\begin{proof}
Note that $\langle z, z \rangle = 0$ implies $z_1^2 = w z_2^2$,
which by our choice of $w$ can only happen when $z = 0$. Thus
Proposition \ref{prop:inner} is satisfied. The group has order
$q^5$, and the three subgroups have size $q^2$, leading to a
pseudo-exponent bound of $2.5$ as claimed.
\end{proof}

A slight variant of this construction works for even $q$ as well,
but the pseudo-exponent bound is identical so we omit the details.

One quite different example is the following Frobenius group of
order $80$.  We found the group by a brute force search using GAP,
and Michael Aschbacher supplied the following humanly
understandable proof that it works.

Let $C_5 \subset \F_{16}^\times$ be the unique subgroup of order
$5$. Consider its semidirect product $G = C_5 \ltimes \F_{16}$
with the additive group of $\F_{16}$, where multiplication is
defined by
$$
(\alpha,x)(\beta,y) = (\alpha\beta, \beta x + y).
$$

\begin{proposition}
\label{proposition:order80} The group $G = C_5 \ltimes \F_{16}$
realizes $\langle 5, 5, 8\rangle$, and thus $\alpha(G) \le
3\log_{200}80 = 2.4811\dots$.
\end{proposition}

\begin{proof}
Let
$$
H_1 = \{(\alpha,0) : \alpha \in C_5\}
$$
and
$$
H_2 = \{(\alpha, \alpha-1) : \alpha \in C_5\}
$$
(i.e., $H_2$ is $H_1$ conjugated by $(1,1)$). Let
$$
H_3 = \{(1,x) : x \in \F_{16}, \Tr x = 0 \},
$$
where $\textup{Tr}$ denotes the trace from $\F_{16}$ to $\F_2$.
These groups satisfy $|H_1| = |H_2| = 5$ and $|H_3| = 8$.  All we
need to check is the triple product property.

We must verify that unless $\alpha$ and $\beta$ are both $1$, the
product
$$
(\alpha,0)(\beta, \beta-1) = (\alpha\beta, \beta-1)
$$
is not in $H_3$.  For it to be in $H_3$, we must have $\alpha =
\beta^{-1}$ and $\Tr (\beta-1) = 0$.  However,
$$
\Tr (\beta - 1) = \Tr \beta - \Tr 1 = \Tr \beta,
$$
and $\Tr \beta = 1$ for $\beta \in C\setminus\{1\}$ because the
minimal polynomial over $\F_2$ of such a $\beta$ is
$1+\beta+\beta^2+\beta^3+\beta^4$.
\end{proof}

This proposition generalizes as follows (see \cite{Brown} for
background on cohomology): Let $G$ be a group that acts on an
abelian group $A$, $\theta : G \to A$ a $1$-cocycle, and $B
\subseteq A$ a subgroup. If $\theta(g) \in B$ implies $g=1$ for
all $g \in G$, then the semidirect product $G \ltimes A$ realizes
$\langle |G|,|G|,|B|\rangle$ via the subgroups $G \times \{0\}$,
$\{(g,\theta(g)) : g \in G\}$, and $\{1\} \times B$.  (In
Proposition~\ref{proposition:order80}, the $1$-cocycle is a
coboundary.) Unfortunately, we do not know any other good
examples.

Unlike the cases of extraspecial groups and matrix groups, we do
not know how to generalize Proposition~\ref{proposition:order80}
to achieve Lie pseudo-exponent arbitrarily near $2$. The best we
know how to do is the following. Let $\Quat$ be the quaternions,
and $U \subset \Quat^\times$ be the group of unit quaternions
(which is isomorphic to $SU(2)$). Then within the semidirect
product $U \ltimes \Quat$, the three subgroups $U \times \{0\}$,
$\{(u,u-1) : u \in U\}$, and $\{(0,x): \Tr x = 0 \}$ satisfy the
triple product property and prove that the Lie pseudo-exponent of
$U \ltimes \Quat$ is at most $7/3$.

\SubSection{Wreath products} \label{section:wreath}

In this section we present another family of groups that achieves
pseudo-exponent $2 + o(1)$. This family is described in terms of
the wreath product: if $A$ is a group, then the wreath product $A
\wr S_n$ is the semidirect product $S_n \ltimes A^n$, where $S_n$
acts on $A^n$ by permuting the coordinates (and the
multiplication is of course via $(\pi, u)(\pi',v) = (\pi\pi',
\pi'u + v)$).

\begin{theorem}
Let $A$ be the cyclic group of order $2n$, and let $G_n = A\wr
S_n$. Then
$$
\alpha(G_n) \le \gamma(G_n) = 2 + \frac{1+\log 2}{\log n} +
O\left(\frac{1}{(\log n)^2}\right).
$$
\end{theorem}

\begin{proof}
We view $G_n$ as the semidirect product $S_n \ltimes A^n$, and
will use the three subgroups
\begin{eqnarray*}
H_1 & = & \{(\pi, 0) : \pi \in S_n\}, \\
H_2 & = & \{(\pi, \pi u - u) : \pi \in S_n\}, \quad\textup{and} \\
H_3 & = & \{(\pi, \pi v - v) : \pi \in S_n\},
\end{eqnarray*}
where $u = (1, 2, \dots, n)$, and $v = (n, n-1, \dots, 1)$.

As each subgroup has size $n!$ in a group of size $n!(2n)^n$,
$$
\alpha \le \frac{\log (n!(2n)^n)}{\log n!},
$$
assuming the triple product property holds. The largest character
degree of $G_n$ is $|S_n| = n!$ (see Theorem~25.6 in \cite{H})
and so $|G|^{1/\gamma} = n!$, which implies
$$
\gamma = \frac{\log (n!(2n)^n)}{\log n!}.
$$
By Stirling's formula,
$$
\frac{\log (n!(2n)^n)}{\log n!} = 2 + \frac{1+\log 2}{\log n} +
O\left(\frac{1}{(\log n)^2}\right),
$$
so all that remains is to verify the triple product property.

Suppose $h_1 = (\pi',0) \in H_1$ and $h_2 = (\pi, \pi u -u ) \in
H_2$.  Their product is $(\pi'\pi,\pi u - u)$, and if it equals
$h_3 = (\sigma, \sigma v - v) \in H_3$, then $\pi u - u = \sigma
v - v$. The $i$-th coordinate of $\pi u - u$ is $\pi(i) - i$, and
that of $\sigma v - v$ is $(n+1 - \sigma(i)) - (n+1-i) =
i-\sigma(i)$. Thus, $h_1h_2=h_3$ implies $\pi(i)+\sigma(i) = 2i$
for all $i$. This is an equation in $A$, and hence holds only
modulo $2n$. However, $\pi(i)$, $\sigma(i)$, and $i$ are all in
$\{1,\dots,n\}$, so the equation holds in the integers as well.
Because $\pi(1)$ and $\sigma(1)$ are both at least $1$, we
conclude from $\pi(1)+\sigma(1)=2$ that $\pi(1)=\sigma(1)=1$.
Then $\pi(2)$ and $\sigma(2)$ must be at least $2$, and
$\pi(2)+\sigma(2)=4$, so $\pi(2)=\sigma(2)=2$, etc.  We conclude
that $\pi$ and $\sigma$ are both trivial, as is $\pi'$ because
$\pi'\pi=\sigma$.  Thus, $h_1=h_2=h_3=1$, as desired.
\end{proof}

This construction is an improvement over
Theorem~\ref{theorem:pseudoexponent2}, because it achieves
essentially the same pseudo-exponent bound, while at the same time
$\alpha \le \gamma$. A more complicated variant of this
construction achieves a comparable pseudo-exponent and has
$\alpha < \gamma$.

\SubSection{Direct products and the Sperner capacity}
\label{section:direct}

It is natural to attempt to improve the pseudo-exponent of a
finite group $G$ by forming some group derived from it, such as a
power $G^k$.  We know that $\gamma(G^k) = \gamma(G)$, so that
parameter becomes no smaller.  Lemma~\ref{lemma:upperbound}
implies that $\alpha(G^k) \le \alpha(G)$, and in this section we
show that it is possible to achieve $\alpha(G^k) < \alpha(G)$.

We will be led for the first time since
Lemma~\ref{lemma:shortexact} to realize matrix multiplication
through quotient sets that are {\em not} subgroups.
Proposition~\ref{prop:dihedral} below proves that this
complication is necessary to determine the pseudo-exponents of
certain groups.

Let $D_m$ be the dihedral group generated by $x$ and $y$, with the
relations $y^2 = x^m = 1$ and $yxy = x^{-1}$.

\begin{proposition}
\label{prop:dihedral} For every $m$, $D_m$ realizes $\langle
2,2,2\lfloor m/3 \rfloor\rangle$, and hence $\alpha(D_m)<3$ for
$m\ge 9$.  If $m$ is a prime greater than $3$, then no three
subgroups prove $\alpha(D_m)<3$.
\end{proposition}

\begin{proof}
Let $S_1 = \langle y \rangle$ be the subgroup generated by $y$,
$S_2 = \langle yx^2 \rangle$, and $S_3 = \{x^{3k}, yx^{3k+1} : 0
\le k < (m-2)/3\}$.  Then one can check by simple case analysis
that $D_m$ realizes $\langle 2,2,2\lfloor m/3 \rfloor\rangle$
through $S_1,S_2,S_3$.  Note that $S_3$ is a subgroup iff $m$ is a
multiple of $3$.

When $m$ is prime, all subgroups of $D_m$ have order $1$, $2$,
$m$, or $2m$, and it is easy to rule out each case (except when
$m=3$, in which case three subgroups of order~$2$ prove
$\alpha(D_3)<3$).
\end{proof}

Proposition~\ref{prop:dihedral} is not optimal: $D_5$ realizes
$\langle 2, 2, 3 \rangle$ through $\{1, y\}, \{1, yx\}, \{1, x^2,
yx^4\}$.  However, we have checked using GAP that it is optimal
for $m=4$, and thus $\alpha(D_4)=3$.

We now use the combinatorial notion of {\em Sperner capacity} to
show that $\alpha(D_4^k)<3$ for large $k$, despite the fact that
$\alpha(D_4)=3$.

\begin{proposition}
If $S \subseteq (\Z/m\Z)^k$ is a subset in which no two distinct
vectors differ by an element of $\{0,1\}^k$, then $D_m^k$ realizes
$\langle 2^k, 2^k, |S|\rangle$.
\end{proposition}

\begin{proof}
We identify $\Z/m\Z$ with the subgroup $\langle x \rangle
\subseteq D_m$ (via $i \leftrightarrow x^i$), so that $S
\subseteq \langle x \rangle^k \subseteq D_m^k$. The subgroups
$\langle y \rangle$ and $\langle yx \rangle$ of $D_m$ have
pointwise product $\langle y \rangle\langle yx \rangle =
\{1,y,yx,x\}$.  Therefore the condition on differences of elements
in $S$ implies that $\langle y \rangle^k$, $\langle yx
\rangle^k$, and $S$ satisfy the triple product property, since
$(\langle y \rangle^k \langle yx \rangle^k) \cap \langle x
\rangle^k = \{1,x\}^k$, and $Q(S) \subseteq \langle x \rangle^k$
avoids $\{1,x\}^k$.
\end{proof}

The problem of making $S$ as large as possible has been studied
before;  a generalization of this problem is known as the Sperner
capacity of a directed graph \cite{Garg,Korn}.  It is known that
$|S| \le (m-1)^k$ (see Theorem~1.2 in \cite{Alon}, which extends
several earlier papers \cite{Blok, Cald}), and that
$$
|S| = (m-1)^{(1-o(1))k}
$$
can be achieved by the following construction:

Assume that $(m-1)$ divides $k$, and take $S$ to be the set of all
vectors in $(\Z/m\Z)^k$ with exactly $k/(m-1)$ occurrences of each
element of $\{0,1,\dots,m-2\}$. Now suppose we have $u, v \in S$
with $u - v \in \{0,1\}^k$.  For each coordinate $i$ such that
$u_i=0$, we have $v_i \in \{0,m-1\}$ because $u_i-v_i \in
\{0,1\}$, and thus $v_i=0$.  Then whenever $u_i=1$, it follows
that $v_i=1$, because all $k/(m-1)$ cases in which $v_i=0$ have
$u_i=0$ as well. Repeating this argument yields $u=v$, as desired.

We conclude that direct products {\em can} help:

\begin{corollary}
We have $\alpha(D_4^k) \le (3 + o(1)) \log_{12}{8}$, which
approaches $3\log_{12}{8} = 2.51\dots$ as $k \rightarrow \infty$.
\end{corollary}

This pseudo-exponent bound comes tantalizingly close to settling
Question~\ref{fundamentalq}: if
$$
\liminf_{k \to \infty} \alpha(D_4^k) < 3 \log_{12} 8,
$$
then the answer to the question is ``yes,'' and our methods do in
fact prove $\omega<3$ (at least).  The same holds in general for
$D_{2n}$ (which has $n-1$ characters of degree $2$ and $4$ of
degree $1$) ; the Sperner capacity construction proves that
$\alpha(D_{2n}^k) \le (3+o(1)) \log_{8n-4} 4n$, and if
$$
\liminf_{k \to \infty} \alpha(D_{2n}^k) < 3\log_{8n-4} 4n,
$$
then the answer to Question~\ref{fundamentalq} is ``yes.''

Also, note that Lemma~\ref{lemma:upperbound} implies that for all
$G$,
$$
\liminf_{k \to \infty} \alpha(G^k) = \inf_{k \ge 1} \alpha(G^k).
$$
Thus, even if the answer to Question~\ref{fundamentalq} is ``no,''
there are combinatorial consequences.  For example, knowing that
$\alpha(D_{2n}^k) \ge 3\log_{8n-4} 4n$ for all $n$ and $k$ would
give a new proof of the Sperner capacity bound $|S| \le (m-1)^k$
above, in the case of even $m$.

\Section{Concluding comments} \label{section:conclusions}

\SubSection{Open questions}

The most pressing question arising in this paper is
Question~\ref{fundamentalq}, which represents a potential barrier
to obtaining non-trivial bounds on $\omega$ using our
techniques.  However, there are numerous other open questions that
are relevant to Question~\ref{fundamentalq} and the ultimate goal
of proving $\omega=2$.

\paragraph{Matrix groups.}
As pointed out in Section~\ref{section:linear}, matrix groups
seem to be one of the most promising families of examples, but we
still know very little about them.  Can our bounds for
$\alpha(SL_2(\F_q))$ be improved?  We see no reason why they
should be optimal.  Recall that beating $9/4$ asymptotically would
settle Question~\ref{fundamentalq}.  We know even less about
$SL_n(\F_q)$ (only that $\alpha(SL_n(\F_q))<3$), so any
non-trivial construction would be of interest.  The only other
finite matrix groups that we have studied are those closely
connected to $SL_n$ (such as $PSL_n$ or $GL_n$), but there are a
number of other families.  What can one say about the
pseudo-exponents of the groups in these families?

\paragraph{Quotient sets.}
The examples in Subsection~\ref{section:direct} show that
quotient sets sometimes outperform subgroups.  For which groups
does this occur?  Are there general constructions of useful
quotient sets other than via Sperner capacity?  Can they be used
to improve our constructions for $S_n$ or the wreath product?
What about matrix groups?

\paragraph{Lie groups.}
Can one use Lie groups to prove anything about $\omega$
directly?  Do results on the Lie pseudo-exponent imply anything
about the pseudo-exponents of related finite groups?  Compact Lie
groups seem more closely analogous to finite groups than
non-compact Lie groups are, so studying them might be
illuminating.  (All of the Lie groups in this paper are
non-compact.)

\paragraph{Group extensions.}
Extensions of groups with pseudo-exponent $3$ can have
substantially smaller pseudo-exponents, as demonstrated by the
solvable groups in Subsection~\ref{section:solvable}.  (Recall
that solvable groups are formed from abelian groups by taking
repeated extensions.)  Is there a general way to lower $\alpha$ or
raise $\gamma$ by taking extensions?  As a first step, can one
find a family of solvable groups with pseudo-exponents tending to
$2$?

\paragraph{Powers of groups.}
The simplest case of group extensions is taking powers of a
group.  Given $G$, what can one say about the \textit{asymptotic
pseudo-exponent\/} $\inf_{k \ge 1} \alpha(G^k)$ of $G$?  As noted
in Subsection~\ref{section:direct}, $\gamma(G^k)=\gamma(G)$, so
if there exists a group such that $\inf_{k \ge 1} \alpha(G^k) =
2$, then $\omega=2$ by Corollary~\ref{cor:race}.

\SubSection{Extensions} \label{section:extensions}

It is natural to attempt to extend our methods in various ways.
For example, one might try to obtain bounds on border ranks of
tensors, perhaps by using deformations of group algebras. It is
also reasonable to ask whether our approach (given its reliance
on representation theory) works in finite characteristic, as well
as over $\C$. As Theorem~\ref{theorem:reduction} indicates, one
can just as easily embed matrix multiplication into $F[G]$ rather
than $\C[G]$, where $F$ has characteristic $p$. As long as $p$
does not divide $|G|$, the representation theory of $G$, and all
other aspects of our approach, work out identically, assuming $F$
is algebraically closed. Sch\"onhage has shown that the exponent
of matrix multiplication over arbitrary fields depends only on the
characteristic (see Corollary~15.18 in \cite{BCS}), so we lose
nothing by requiring that $F$ be algebraically closed.

We conclude by mentioning a particular variant of our approach
that does not require any control of the character degrees, and
thus may still be viable even if there is a negative answer to
Question~\ref{fundamentalq}. We have found less structure to make
use of, and it seems less attractive, but it uses similar ideas.
Suppose we have distinct elements $x_{i, j}, y_{k, \ell} \in G$,
for $1\le i \le n$, $1\le j,k \le m$, and $1\le \ell \le p$, such
that
\begin{equation}
x_{i,j}y_{j, \ell} \sim x_{i', k}y_{k', \ell'} \; \Leftrightarrow
\; i=i', k=k', \ell = \ell', \label{eq:conjugacy}
\end{equation}
where $\sim$ denotes conjugacy of elements. Then we embed matrix
$A = (a_{i, j})$ as $\bar{A} = \sum_{i, j} a_{i, j}x_{i, j} \in
\C[G]$, and matrix $B = (b_{k, \ell})$ as $\bar{B} = \sum_{k,
\ell} b_{k, \ell}y_{k, \ell} \in \C[G]$.  We can pursue a similar
strategy to compute $AB$. In this case, however, in the Fourier
domain, we need only to compute the {\em trace} of each of the
matrix products in the block-diagonal matrix multiplication. That
requires only $\sum_i d_i^2 = |G|$ multiplications, and so we can
conclude that the rank of $\langle n,m,p \rangle$ is at most
$|G|$.

Let $G$ be a group with subsets $S_1, S_2$ and $S_3$ satisfying
the triple product property. If we replace {\em conjugacy} with
{\em equality} in \eqref{eq:conjugacy}, then it can be satisfied
by taking $\{x_{i, j}\} = S_1 S_2^{-1}$ (where $i$ indexes $S_1$
and $j$ indexes $S_2$) and $\{y_{k, \ell}\} = S_2 S_3^{-1}$ ($k$
indexes $S_2$ and $\ell$ indexes $S_3$), so it is possible that
the techniques we have developed in this paper could help with
this variant as well, although in general we find it difficult to
work with conjugacy constraints.

\section*{Acknowledgements}

We are grateful to Michael Aschbacher, Noam Elkies, Bobby
Kleinberg, L\'aszl\'o Lov\'asz, Amin Shokrollahi, David Vogan, and
Avi Wigderson for helpful discussions.


\begin{thebibliography}{10}\setlength{\itemsep}{-1ex}\small

\bibitem{Alon}
N.~Alon.
\newblock On the capacity of digraphs.
\newblock {\em European J.\ Combinatorics}, 19:1--5, 1998.

\bibitem{Blok}
A.~Blokhuis.
\newblock On the {Sperner} capacity of the cyclic triangle.
\newblock {\em J.\ Algebraic Combinatorics}, 2:123--124, 1993.

\bibitem{Brown}
K.~S. Brown.
\newblock {\em Cohomology of Groups}.
\newblock Number~87 in Graduate Texts in Mathematics. Springer-Verlag, 1982.

\bibitem{BCS}
P.~{B\"urgisser}, M.~Clausen, and M.~A. Shokrollahi.
\newblock {\em Algebraic Complexity Theory}, volume 315 of {\em Grundlehren der
  mathematischen Wissenschaften}.
\newblock Springer-Verlag, 1997.

\bibitem{Cald}
R.~Calderbank, P.~Frankl, R.~L. Graham, W.~Li, and L.~Shepp.
\newblock The {Sperner} capacity of the cyclic triangle for linear and
  nonlinear codes.
\newblock {\em J.\ Algebraic Combinatorics}, 2:31--48, 1993.

\bibitem{CW}
D.~Coppersmith and S.~Winograd.
\newblock Matrix multiplication via arithmetic progressions.
\newblock {\em J. Symbolic Computation}, 9:251--280, 1990.

\bibitem{GAP}
The GAP~Group.
\newblock {\em {GAP -- Groups, Algorithms, and Programming, Version 4.3}},
  2002.
\newblock \texttt{(http://www.gap-\break system.org)}.

\bibitem{Garg}
L.~Gargano, J.~{K\"orner}, and U.~Vaccaro.
\newblock {Sperner} theorems on directed graphs and qualitative independence.
\newblock {\em J.\ Combinatorial Theory Series A}, 61:173--192, 1992.

\bibitem{Green}
J.~A. Green.
\newblock The characters of the finite general linear groups.
\newblock {\em Transactions of the American Mathematical Society}, 80:402--447,
  1955.

\bibitem{H}
B.~Huppert.
\newblock {\em Character Theory of Finite Groups}.
\newblock Number~25 in de Gruyter Expositions in Mathematics. Walter de
  Gruyter, Berlin, 1998.

\bibitem{JL}
G.~James and M.~Liebeck.
\newblock {\em Representations and Characters of Groups}.
\newblock Cambridge University Press, Cambridge, second edition, 2001.

\bibitem{Korn}
J.~{K\"orner} and G.~Simonyi.
\newblock A {Sperner}-type theorem and qualitative independence.
\newblock {\em J.\ Combinatorial Theory}, 59:90--103, 1992.

\bibitem{LR}
J.~Lafferty and D.~Rockmore.
\newblock Fast fourier analysis for {$SL_2$} over a finite field and related
  numerical experiments.
\newblock {\em Experimental Mathematics}, 1:115--139, 1992.

\bibitem{S}
V.~Strassen.
\newblock Gaussian elimination is not optimal.
\newblock {\em Numerical Mathematics}, 13:354--356, 1969.

\bibitem{VK}
A.~M. Vershik and S.~V. Kerov.
\newblock Asymptotics of the largest and the typical dimensions of irreducible
  representations of a symmetric group.
\newblock {\em Functional Analysis and its Applications}, 19:21--31, 1985.

\end{thebibliography}
\end{document}